\documentclass[11pt]{article}
\usepackage{latexsym,amssymb,amsmath,amsthm,graphicx,float,calligra,authblk}

\newcommand{\reals}{\mathbb{R}}

\newcommand{\eps}{\varepsilon}

\newcommand{\pd}{\partial}
\newcommand{\vf}{\varphi}

\newcommand{\<}{\langle}
\renewcommand{\>}{\rangle}

\newcommand{\tr}{\mbox{Tr}}
\newcommand{\E}{\mathbb{E}}

\newtheorem{theorem}{Theorem}
\newtheorem{lemma}{Lemma}

\title{Matrix probing: \\ a randomized preconditioner for the wave-equation Hessian}
\author[1]{Laurent Demanet}
\author[2]{Pierre-David L\'{e}tourneau}
\author[3]{\\ Nicolas Boumal}
\author[4]{Henri Calandra}
\author[1]{Jiawei Chiu}
\author[5]{Stanley Snelson}
\affil[1]{Department of Mathematics, MIT}
\affil[2]{Institute for Computational Mathematics and Engineering, Stanford}
\affil[3]{Applied Mathematics Department, Universit\'{e} catholique de Louvain}
\affil[4]{Total SA, Exploration \& Production}
\affil[5]{Courant Institute of Mathematical Sciences, NYU}

\date{December 2010}		

\begin{document}
\maketitle

\begin{abstract}
This paper considers the problem of approximating the inverse of the wave-equation Hessian, also called normal operator, in seismology and other types of wave-based imaging. An expansion scheme for the pseudodifferential symbol of the inverse Hessian is set up. The coefficients in this expansion are found via least-squares fitting from a certain number of applications of the normal operator on adequate randomized trial functions built in curvelet space. It is found that the number of parameters that can be fitted increases with the amount of information present in the trial functions, with high probability. Once an approximate inverse Hessian is available, application to an image of the model can be done in very low complexity. Numerical experiments show that randomized operator fitting offers a compelling preconditioner for the linearized seismic inversion problem.
\end{abstract}

{\bf Acknowledgments}. LD would like to thank Rami Nammour and William Symes for introducing him to their work. LD, PDL, and NB are supported by a grant from Total SA.

\section{Introduction}

\subsection{Problem setup: Gauss-Newton iterations}

This paper considers the imaging problem of determining physical characteristics in a region of space given surface measurements of scattered waves. Several imaging modalities fall under this umbrella (ground-penetrating radar, nondestructive acoustic testing, remote personnel assessment), but in the sequel we focus exclusively on the example of reflection seismology. Throughout this paper we let $m(x)$ for the physical parameters in the subsurface, and $d(r,s,t)$ for the recorded waveforms (seismograms). Here $x$ are the space coordinates in the volume, $r$ is receiver position, $s$ is source position, and $t$ is time.

A most popular way of treating the inversion problem of recovering $m$ from $d$
is through the minimization of the output least-squares functional
\[
J[m] = \frac{1}{2} \| d - \mathcal{F}[m] \|^2_2,
\]
where $\mathcal{F}$ is the nonlinear map for predicting data from the model $m$. In this paper we restrict ourselves to the setup of constant-density acoustics, and let $m(x)$ be a variable wave speed. The prediction $\mathcal{F}[m]$ then consists of the solutions $u$ -- sampled at the receivers $(r,t)$ -- of the acoustic wave equations
\[
m(x) u_{tt} - \Delta_x u = f_s,
\]
with different right-hand sides $f_s(x,t)$ index by $s$ (the source). The notation $\| \cdot \|^2_2$ refers to the sum of the squares of the components. The quantity $m(x)$ is the inverse of the square of the local wave speed.

Whether data is considered all at once, or by frequency increments as in full waveform inversion, the procedure for minimizing $J[m]$ is usually some variant of the Gauss-Newton method, which consists in linearizing $J[m]$ about some current vector $m_0$. Specifically, if a new vector $m_1$ is sought so that $J[m_1]$ is closer to the minimum than $J[m_0]$ is, then we first write
\[
J[m_1] = J[m_0] + \< \frac{\delta J}{\delta m}[m_0] , \, \delta m \> + \frac{1}{2} \< \delta m , \,  \frac{\delta^2 J}{\delta m^2}[m_0] \, \delta m \> + ...
\]
where $\delta m = m_1 - m_0$, and find $\delta m$ as the minimum of the quadratic form above. The solution is
\[
0 = \frac{\delta J}{\delta m}[m_0] + \frac{\delta^2 J}{\delta m^2}[m_0] \, \delta m \qquad \Rightarrow \qquad \delta m = - \left( \frac{\delta^2 J}{\delta m^2}[m_0] \right)^{-1} \frac{\delta J}{\delta m}[m_0].
\]
This equation is a Newton descent step: it is then applied iteratively to obtain a new $m_2$ from $m_1$, etc.  The Hessian is the operator $\frac{\delta^2 J}{\delta m^2}[m_0]$.

\bigskip

If $J[m]$ is the least-squares misfit functional above, then by denoting
\[
F = \frac{\delta \mathcal{F}}{\delta m}[m_0],
\]
we obtain the first and second variations of $J$ as
\[
\frac{\delta J}{\delta m}[m_0] = - F^* ( d - \mathcal{F}[m_0] ),
\]
\[
\frac{\delta^2 J}{\delta m^2}[m_0] = F^* F - \< \frac{\delta^2 \mathcal{F}}{\delta m^2}[m_0], d - \mathcal{F}[m_0] \>.
\]
The migration operator $F^*$ acts from data space to model space, and is most accurately computed by reverse-time migration. The demigration operator $F$ acts from model space to data space, and can be computed by solving a forward ``modeling" wave equation. The term involving the second variation of $\mathcal{F}$ in the expression of the Hessian is routinely discarded on the basis that $\mathcal{F}$ is ``locally well-linearized" -- a heuristically plausible claim when $m_0$ is smooth in comparison to $\delta m$ -- but which has so far eluded rigorous analysis. With this simplification in mind, we refer to the (reduced) Hessian $H$ as the leading-order contribution
\[
H = F^* F.
\]
This linear operator is also called the normal operator, and acts within model space. The Newton descent step then calls for computing the pseudoinverse $F^{+} = (F^* F)^{-1} F^*$, well-known to arise in the solution of the overdetermined linearized least-squares problem.

\bigskip

Physically, inversion of the Hessian corresponds to the idea of correcting for low levels of illumination of the medium by the forward (physical) wavefield. Although illumination seems to make good sense as a function of space $x$, it is in fact unclear how to define it as such. Rather, it is more appropriate to define illumination as a function in \emph{phase-space}, i.e., the set of $x$ and $k$ (wave vectors). In the words of Nammour and Symes \cite{Nammour2}, illumination is not just a scaling, but a dip-dependent scaling. This paper follows this idea by considering the pseudodifferential symbol of the Hessian.

\bigskip

While reasonably efficient methods of applying the operators $F$ and $F^*$ to vectors are common knowledge, little is currently known about the structure of the inverse Hessian $H^{-1}$. Direct linear algebra methods for computing a matrix inverse are out of the question, because the matrix $H$ is too large to be formed in practice. This also prevents the immediate application of methods such as BFGS. 
Iterative linear algebra methods such as GMRES or LSQR can be set up, but need a very large number of iterations to converge due to the poor conditioning of $H$.  The problem of slow convergence is particularly acute since the full prestack data space (after the application of $F$) is much larger than poststack model space -- hence each application of $H = F^*F$ is very costly. The obvious alternative to the Gauss-Newton iteration, namely straight gradient descent without considering the Hessian, is even less attractive than GMRES for solving the ill-conditioned linearized least-squares problem. 

\bigskip

Preconditioning is needed to properly guide the inversion iterations. A preconditioner for a matrix $H$ is a matrix $M$ that approximates the inverse $H^{-1}$. It can be used to rewrite $H x = b$ as
\[
PHx = PB,
\]
where now only the matrix $PH$ needs to be inverted. An alternative formulation where $P$ postmultiplies $H$ is also possible. Several preconditioners for the wave equation Hessian have already been proposed in the literature: they are reviewed in context in Section \ref{sec:review}.

\bigskip

This paper solves the preconditioning problem by ``probing", or testing the Hessian by applying it to a small number of randomized vectors, followed by a fit of the inverse Hessian in a special expansion scheme in phase-space. Our work is closest in spirit to that of Nammour and Symes \cite{Nammour1, Nammour2} and Herrmann et al. \cite{HerrmannBrown} (which in turn follows from a legacy of so-called scaling preconditioners reviewed below) but departs from it in that the trial space is randomized instead of being the Krylov subspace of the migrated model\footnote{The Krylov subspace of a vector $y$, for a matrix $H$, is the space spanned by $y, Hy, H^2 y$, etc.}. Randomness of the trial functions guarantees recovery of the action of the inverse Hessian on a much larger linear subspace than is normally the case with a deterministic method. This claim is backed both by numerical experiments (Section \ref{sec:num}) and by a theoretical justification (Section \ref{sec:theory}).

\bigskip

The proposed approach bridges a gap in the literature, in that we obtain quantitative results -- hence finally a rationale -- for the probing methods to precondition the wave-equation Hessian. We found that randomization is an important step to achieve such guarantees, and may be an attractive numerical choice in its own right.

\subsection{Pseudodifferential symbol of the Hessian}\label{sec:psido}

To provide an expansion scheme for the inverse Hessian, it is important to understand its structure as a pseudodifferential operator. In the sequel we consider only two spatial dimensions $x = (x',z)$, but the main ideas do not depend on this assumption. 

\bigskip

It is well-known that migration $F^*$ is a ``kinematic" inverse of the modeling operator $F$ in the sense that the mapping of singularities generated by $F^*$ generically undoes that of $F$. Putting technical pathologies aside, this claim means that $H = F^* F$ does not change the location of singularities in model space. Hence the Hessian is ``microlocally equivalent" to the identity, or ``microlocal" for short. This property was understood and made precise by at least the following people.
\begin{itemize}
\item In 1985, Beylkin showed that the Hessian is pseudodifferential in the absence of caustics, and in the context of generalized Radon transforms \cite{Beylkin}.
\item In 1988, Rakesh removed the no-caustic assumption, but considers a point source and full-aperture (whole-Earth) data \cite{Rakesh} .
\item In 1998, ten Kroode, Smit and Verdel showed that Beylkin's result still holds if a less restrictive ``traveltime injectivity condition" is satisfied \cite{tenKroode} .
\item In 2000, Stolk refined these results by showing that the Hessian is generically invertible: if a $C^\infty$ wave speed does not give rise to a pseudodifferential Hessian, an arbitrarily small $C^\infty$ perturbation of it will \cite{Stolk}.
\end{itemize}

The consequence of this body of theory for the problem of designing a compressed numerical representation of the Hessian is the following. We will consider a representation of the Hessian as a pseudodifferential operator:
\begin{equation}\label{eq:symbol}
H m(x) = \int e^{ix \cdot k} a(x,k) \, \hat{m}(k) \, dk,
\end{equation}
where hat denotes Fourier transformation in the spatial variables.  The amplitude, or symbol $a(x,k)$ plays the role of illumination in phase-space $(x,k)$ as alluded to earlier.

\bigskip

There is nothing special about writing an integral sign instead of a sum: interpolation and sampling allow to transform number arrays into functions and vice-versa. Keeping $x$ and $k$ continuous for the time being however offers the opportunity to discuss the important point: \emph{smoothness} of the symbol $a(x,k)$. Indeed, while the symbol representation (\ref{eq:symbol}) is always available whichever the linear operator considered, the symbol will be smooth in a very specific way for ``microlocal" operators as discussed above. We say that the symbol $a(x,k)$ is of order $e$ (and type $(1,0)$) if it obeys the condition
\begin{equation}\label{eq:type10}
| \pd_k^\alpha \pd_x^\beta a(x,k) | \leq C_{\alpha\beta} (1 + |k|^2)^{(e - |\alpha|)/2},
\end{equation}
where $\alpha = (\alpha_1, \alpha_2)$, $\pd_k^\alpha = \frac{\pd^{\alpha_1}}{\pd k_1^{\alpha_1}} \frac{\pd^{\alpha_2}}{\pd k_2^{\alpha_2}}$, $|\alpha| = \alpha_1 + \alpha_2$, and similarly for $\beta$. Notice that since the bound \emph{decreases} by one power of $(1 + |k|^2)^{1/2}$ for every derivative in $k$, it means that the larger $|k|$ the smoother the symbol $a(x,k)$. If we had considered the symbol of either $F$ or $F^*$ instead, each derivative in $k$  space would have \emph{increased} the value of the symbol by a quantity proportional to $k$. Physically, illumination is a phase-space concept, but it is ``not too far" from being purely a function of $x$ since the $k$ dependence of $a$ is extremely smooth for large $|k|$.

\bigskip

There is one very idealized scenario in which the Hessian obeys the condition (\ref{eq:type10}) with order $e = 1$. The assumptions are the following: 1) sufficiently fine Cartesian sampling of the data in time and receiver coordinate (so that the sum can easily be written as an integral), 2) full aperture acquisition, 3) a point-impulse wavelet $\delta(t)$ in time, and 4) smooth and generic\footnote{As above, ``generic" here refers to the absence of kinematic exceptions that would discredit migration as a microlocal inverse, as discussed in \cite{Stolk}. Smooth means infinitely differentiable with oscillations on a length scale much larger than the wavelength of the wave. Random smooth media are ``generic" with probability 1.} background physical parameters such as wave speed. If all these conditions are met, it is known at least from \cite{Stolk} that the Hessian has a symbol that obeys (\ref{eq:type10}).

\bigskip

In turn, if a symbol obeys (\ref{eq:type10}), it is by now well-known that it is in fact extraordinarily compressible numerically. Bao and Symes \cite{Bao} show that the asymptotic behavior of $a(x,k)$ as $k \to \infty$, i.e. the action of the Hessian at very small scales, can be encoded using only a few Fourier series coefficients in $x$ and in $\theta = \arg k$:
\begin{equation}\label{eq:compr1}
a(x,k) \sim \sum_{\lambda, q} c_{\lambda, n} \, e^{i \lambda \cdot x} \, e^{iq\theta} \, |k|
\end{equation}
Recent work by Demanet and Ying \cite{DSC} has shown how to add degrees of freedom in the radial wave number variable $|k|$ to obtain an $\epsilon$-accurate expansion of $a(x,k)$:
\begin{equation}\label{eq:compr2}
a(x,k) = \sum_{\lambda, q_1, q_2} c_{\lambda, q_1, q_2} \, e^{i \lambda \cdot x} \, e^{i q_1 \theta} \, TL_{q_2}(|k|) \, |k|  + O(\epsilon),
\end{equation}
where the $TL$ are rational Chebyshev functions. The number of terms in the sum is a $O(\epsilon^{-M})$ for all $M>0$. Other expansion schemes exist, such as the hierarchical spline grids in $k$ space, considered in \cite{DSC}. In practice, symbols are considered for values of $k$ that obey $\max \{ |k_1|, |k_2| \} \leq \pi \sqrt{n}$ for some large $n$. In view of the Shannon sampling theorem, this restriction corresponds to sampling (2D) functions on a square grid as vectors of length $n$, and operators such as the Hessian as matrices of size $n$-by-$n$. Both (\ref{eq:compr1}) and (\ref{eq:compr2}) are good approximations of the symbol $a(x,k)$ in the sense that they each contain a number of terms \emph{independent of the size $n$ of the matrix} that eventually realizes the Hessian. 

\bigskip

In three spatial dimensions, spherical harmonics would be used in place of complex exponentials in angle. Otherwise, the symbol expansion scheme needs not be changed.

\bigskip

Equation (\ref{eq:compr2}) provides a decomposition of $H$ into ``elemetary operators" $B_i$, each with symbol $e^{i \lambda \cdot x} \, e^{in\theta} \, TL_q(|k|) \, |k|$. The index $i$ is a shorthand for $(\lambda, q_1, q_2)$, and accordingly we let $b_i$ for the coefficients $c_{\lambda, q_1, q_2}$. In this more compact notation we have the fast-converging expansion
\[
H = \sum_i b_i B_i
\]
for the Hessian.

\bigskip

It is not the Hessian that is of interest, but rather the inverse Hessian. Fortunately, it is a result of Shubin \cite{Shubin} that if the symbol $a(x,k)$ of an operator obeys (\ref{eq:type10}), and if this operator is assumed to be invertible, then the symbol $b(x,k)$ of the inverse of the operator also obeys (\ref{eq:type10}), namely
\begin{equation}\label{eq:type10_b}
| \pd_k^\alpha \pd_x^\beta b(x,k) | \leq D_{\alpha\beta} (1 + |k|^2)^{(-1- |\alpha|)/2},
\end{equation}
with constants $D_{\alpha\beta}$ that are possibly different from $C_{\alpha\beta}$. Notice that the order is now $-1$. In other words, smoothness of the symbol is preserved, or closed, under inversion. If the operator is invertible but only barely so (small singular values which are not regularized), then the constants $D_{\alpha\beta}$ may become large, but the behavior under differentiations in $k$ space is still controlled by (\ref{eq:type10_b}). Note in passing that $b(x,k)$ is not exactly given by $1/a(x,k)$, but the latter is an approximation of $b(x,k)$ that mathematicians find satisfying when $|k|$ is large.

\bigskip

Using the same expansion scheme as above we write 
\[
H^{-1} = \sum_i c_i B_i,
\]
with different coefficients $c_i$.

\subsection{Ill-conditioning}

The four assumptions on the sampling, aperture, wavelet, and medium enumerated earlier are of course far from being realistic in practice. Their violation invariably creates ill-conditioning in the form of a linear subspace in model space where applying the Hessian will return very small values.  This issue manifests itself as small values of the symbol $a(x,k)$. For instance (and this may not be an exhaustive list),
\begin{itemize}
\item Limiting the sampling and the aperture will create an angular deficiency in the sense that reflectors with certain orientations will not be visible in the dataset. The symbol $a(x,k)$ will take on small values for the kinematically ``invisible" $x$ and $k$.
\item Restricting the wavelet in $\omega$ space (frequency) will have the effect to remove low and high wavenumbers from the data. This will have the effect of restricting the symbol $a(x,k)$ in wave number $|k|$.
\item Finally, complicated kinematics of the background wave speed(s) may create shadow zones in which there is very poor illumination. Such is the region behind an impenetrable sphere. In that case the symbol $a(x,k)$ becomes very small in those inaccessible regions.
\end{itemize}

\begin{center}
\begin{figure}[H]
\centering
\includegraphics[width=9cm]{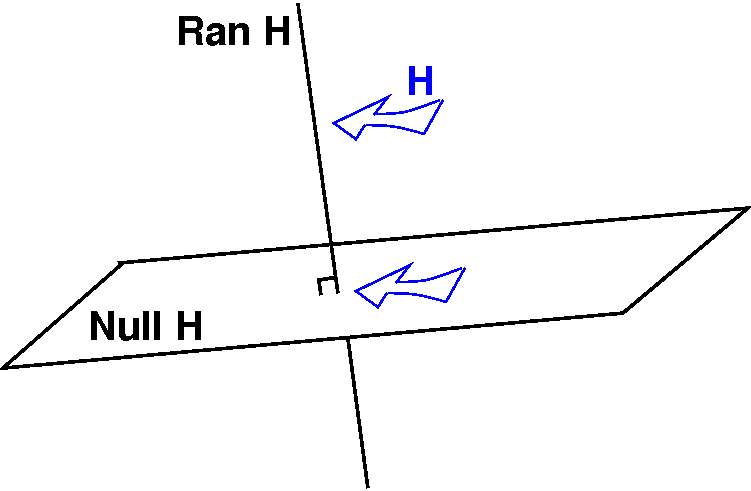}
\caption{The fundamental subspaces of the Hessian (or of any Hermitian matrix). The nullspace is the set of vectors in model space to which an application of the Hessian produces zero values, and whose information is lost. The range space is the set of vectors which are the image of some other vector through an application of the Hessian --- its dimension is the rank of $H$. The blue arrows indicate that, under the action of $H$, the whole space gets mapped to the range space, while the nullspace gets mapped to the origin.}
\label{fig:2spaces}
\end{figure}
\end{center}

The subspace of model space in which the Hessian produces small values is a numerical version of its nullspace.\footnote{The ``numerical nullspace" is precisely defined as the span of the right singular vectors corresponding to singular values below some threshold.} Because this subspace is nonempty, not all vectors in model space are accessible from applying the Hessian to some other vector: the range space does not have full dimension. In other words, the Hessian does not have full rank. \footnote{The ``numerical range space"  is precisely defined as the span of the left singular vectors corresponding to singular values above some threshold.} It is well-known from linear algebra that the dimension of the (numerical) nullspace is equal to the codimension of the (numerical) range space. Because the Hessian is symmetric, the range space is in fact orthogonal to the nullspace -- and ditto of their numerical versions. Figure \ref{fig:2spaces} depicts the fundamental subspaces of the Hessian.

In spite of these complications, this paper speculates that for models well \emph{inside the range space} of the Hessian, an estimate like (\ref{eq:type10_b}) for the inverse Hessian holds. It is not currently known whether this is true theoretically, but we show numerical evidence that supports the claim.

\subsection{Randomized fitting}\label{sec:randfit}

We now address the question of fitting the coefficients in an expansion scheme for the symbol of the inverse Hessian, from application of the Hessian on randomized trial functions. For the time being we assume that the Hessian is invertible and well-conditioned; we return to the discussion of the nullspace in the next section.

\bigskip

Assume that the inverse Hessian is an $n$-by-$n$ matrix that can be expanded as
\begin{equation}\label{eq:expand}
H^{-1} = \sum_{i=1}^p c_i B_i,
\end{equation}
where $B_i$ are themselves matrices, and $p$ counts the number of terms. One possible choice for the $B_i$ was given in Section \ref{sec:psido} (up to discretization), but here the discussion is general. Denote by $\mathbf{y}$ a vector of independent and identically distributed (i.i.d.) Gaussian random variables, in model space -- a ``noise" vector. The application of the Hessian to $\mathbf{y}$ is available:
\[
\mathbf{x} = H \mathbf{y} \qquad \Leftrightarrow \qquad \mathbf{y} = H^{-1} \mathbf{x}.
\]
Given this information, we may now solve for the coefficients $c_i$ in
\[
\mathbf{y}_j = \sum_{i = 1}^p c_i (B_i \mathbf{x})_j, \qquad j = 1, \ldots, n.
\]
This linear system can be overdetermined only if $p \leq n$; in that case the least-squares solution is
\[
c_i = \sum_{j, k} (M^{-1})_{ij} (B_j \mathbf{x})_k \mathbf{y}_k,
\]
where
\[
M_{ij} = \mathbf{x}^T B^T_i B_j \mathbf{x}.
\]
The coefficients $c_i$ can therefore be solved for, in a unique and stable manner, provided the matrix $M$ is invertible and well-conditioned. As we show in the sequel, the invertibility of $M$ hinges on two important assumptions on the elementary matrices $B_i$:
\begin{enumerate}
\item The $B_i$ obey an $H$-dependent near-orthogonality relation:
\[
\E M_{ij} = \tr ( H B_i^T B_j H) \simeq \delta_{ij},
\]
which we express more precisely as requiring that $\E M_{ij}$ be positive definite. The symbol $\E$ stands for mathematical expectation, or ``average over an infinite number of random realizations".
\item Each $B_i$ is a full-rank (invertible), well-conditioned matrix.
\end{enumerate}

When those two conditions are met, we show in section \ref{sec:theory} that $M_{ij}$ is an invertible matrix, with high probability, provided $p$ is large enough, on the order of the square root $\sqrt{r}$ of the rank $r$ of $M$. This result may not be tight but has the advantage of motivating the two assumptions above. We suspect that the number $p$ of coefficients $c_i$ that can be fitted with this method is in fact closer to a constant times $r / \log^2 r$ -- this will be the subject of a separate study. 

\bigskip

The expansion schemes in equations (\ref{eq:compr1}) and (\ref{eq:compr2}) correspond to matrices $B_i$ that obey the above conditions. 

\bigskip

Notice that if the expansion (\ref{eq:expand}) is accurate, i.e. that $H^{-1}$ is determined as a linear combination of the $B_i$, then the proposed method recovers the \emph{whole matrix} $H^{-1}$ in compressed form, not just the action of the matrix $H^{-1}$ on the trial vector $\mathbf{x}$. This property is important: we call it generalizability. The action of $H^{-1}$ can be reliably ``generalized" from its knowledge on $\mathbf{x}$, to other vectors. The randomness of the vector $\mathbf{y}$ is essential in this regard: it would be much harder to argue generalizability if the vector $\mathbf{y}$ had been chosen deterministically. The numerical experiments in section \ref{sec:num} confirm this observation.

\bigskip

Finally, it is worth noting that $H^{-1}$ needs not be given exactly by a sum of $p$ terms of the form $c_i B_i$. If the series converges fast instead of terminating exactly, it is possible to show that the coefficients $c_i$ are determined up to an error commensurate with the truncation error of the series.

\subsection{Fitting via randomized curvelet-based models}\label{sec:curvelettrial}

As mentioned earlier, inversion of the wave-equation Hessian is complicated by various factors that create ill-conditioning. The lack of invertibility not only prevents randomized fitting to work as presented in the previous section, but it also adds to the numerical complexity of the inverse Hessian itself. Just being able to specify the numerical nullspace -- the subspace in which the Hessian erases information -- is at least as complex as specifying the action of the inverse Hessian away from it. As a consequence, it may be advantageous for a coarse preconditioner not to explicitly try and invert the Hessian in the neighborhood of the numerical nullspace.

\bigskip

Our solution to the ill-conditioning problem is to consider noise realizations $y$ that avoid the nullspace, i.e., belong to the range space of $H$. The relation $\mathbf{y} = H^{-1} \mathbf{x}$ then makes sense if we understand $H^{-1}$ as the pseudo-inverse of $H$. The numerical nullspace of $H$ is best described in phase space: it corresponds to the points $(x,k)$ where the symbol $a(x,k)$ of $H$ is small.  This calls for considering an illumination mask, i.e., a simple 0-1 function which indicates whether a point $(x,k)$ is in the essential support of the symbol (value 1) or not (value 0). This piece of a priori information is then used to filter out components of the noise vector (in $x$ space) which would otherwise intersect the nullspace of the Hessian.
 
\bigskip

An explicit expression for the pseudo-differential operator $H$ can be obtained in the idealized case of densely sampled data with idealized sources and receivers. The process involves the asymptotic expansion(stationary phase analysis) of a Generalized Radon Transform and is described in \cite{Beylkin}. We use this expansion as a way of isolating the null space of $H$.

\bigskip

Concretely, we built this illumination indicator function in curvelet-transformed model space.  Curvelets are directional generalizations of wavelets which are efficient at representing bandlimited  wavefronts in a sparse manner \cite{FDCT, 3DFDCT}, and have had applications for regularizing the inversion in seismic imaging \cite{FDCT, HerrmannBrown, HerrmannStolk}. They also provide a sparse representation of wave propagators \cite{CandesDemanet}.  Each curvelet $\vf_\mu(x)$  is indexed by a position vector $x_\mu$ and  a wave vector $k_\mu$. Any (square-integrable) function $f$ can be expanded in curvelets as
\[
f(x) = \sum_\mu f_\mu \vf_\mu(x), \qquad  f_\mu = \int \overline{\vf_\mu}(x) f(x) \, dx.
\]
As explained in Section \ref{sec:theory2}, curvelets efficiently discriminate between different regions of phase-space where the symbol of the Hessian takes on different values.

\bigskip

Consider $S$, the set of curvelets $\vf_\mu$ whose center $(x_\mu, k_\mu)$ belongs to the essential support of the symbol $a(x,k)$ of the Hessian. The stationary phase analysis mentioned above \cite{Beylkin, Symes-notes} reveals the geometric interpretation of these phase-space points: they are \emph{visible}, in the (microlocal) sense that there is a ray linking some source $s$ to the point $x$, reflecting at $x$ in a specular fashion about the normal vector $k$, and then linking $x$ back to some receiver $r$. When a curvelet is visible, it means that it acts like a ``local reflector" for some waves that end up being observed in the dataset. More precisely, a phase-space point $(x_\mu, k_\mu)$ belongs by definition to $S$ if there exist two rays $\gamma_s, \gamma_r$ originating from $x_\mu$ such that:
\begin{itemize}
\item $\gamma_s$ links $x_\mu$ to some source in the source manifold\footnote{Conventionally, an interval or an otherwise open set of positions in which the sampling of sources (resp. receivers) is dense enough in view of the typical wavelength of the seismic waves.};
\item $\gamma_r$ links $x_\mu$ to some receiver in the receiver manifold; and
\item $\gamma_r$ is a reflected ray for $\gamma_s$ at $x_\mu$, i.e., the angle of incidence is equal to the angle of reflection and the two rays form a plane with the normal direction $k_\mu$. 
\end{itemize}
The rays are obtained by ray-tracing from the Hamiltonian system of geometrical optics. The illumination mask is then the sequence equal to $1$ if $\mu \in S$, and zero otherwise. A noise realization $\mathbf{y}$ in curvelet space, filtered by the illumination mask, is simply
\[ 
\mathbf{y}_\mu = \left\{ \begin{array}{ll}
         N(0,\sigma^2 \| \vf_\mu \|_2^2) \; \mbox{i.i.d.} & \mbox{if $\mu \in S$};\\
         0 & \mbox{if $\mu \notin S$}.
         \end{array} \right. 
\] 
The sequence $y_\mu$ is then inverted to yield
\[
\mathbf{y} = \sum_\mu \vf_\mu(x) \mathbf{y}_\mu.
\]
The rest of the algorithm for determining the inverse Hessian then proceeds as in the previous section.

\bigskip

Once the inverse Hessian is available as a series (\ref{eq:expand}), the algorithm for applying it to a vector like the migrated model is well-known and very fast \cite{Bao, DSC}.

\subsection{Previous work}\label{sec:review}

Being able to extract information on the inverse Hessian from a single application of the Hessian is a very good idea which perhaps first appeared, in seismology, in the work of Claerbout and Nichols \cite{Claerbout}. There, a single scalar function of $x$ is sought to represent inverse illumination. In our notations, they seek to fit a symbol $b(x,k)$ which is not a function of $k$.

\bigskip

This work generated refinements that W. Symes puts under the umbrella of  ``scaling methods". In 2003, Rickett \cite{Rickett} offers a solution similar to that of Claerbout and Nichols. In 2004, Guitton \cite{Guitton} proposes a solution based on ``nonstationary convolutions" which corresponds to considering a symbol $b(x,k)$ which is essentially only a function $k$. In 2008, Symes \cite{Symes} proposes to consider symbols of the form
\[
b(x,k) \sim f(x) |k|^{-1},
\]
i.e. which have the proper homogeneity behavior in $|k|$. In 2009, Nammour and Symes \cite{Nammour1, Nammour2} upgrade to the Bao-Symes expansion scheme given in equation (\ref{eq:compr1}). In 2009, Herrmann et al. \cite{HerrmannBrown} propose to realize the scaling as a diagonal operator in curvelet space. 

\bigskip

In all these papers, it is the remigrated image to which the inverse Hessian is applied; in contrast, our paper uses randomized curvelet trial functions. For the representation of the inverse Hessian, we use both (\ref{eq:compr1}) and (\ref{eq:compr2}) for its symbol.

\bigskip

It should also be noted that Herrmann et al. \cite{Herrmann-2003} already proposed in 2003 to realize a curvelet-diagonal approximation of the Hessian, obtained by randomized testing of the Hessian.

\bigskip

The idea of recovering a matrix that has a given sparsity pattern or some other structure from a few applications on well-chosen vectors (``probing") also appeared in the 1990 work of Chan and Keyes on domain-decomposition preconditioning for convection-diffusion problems \cite{Chan1}. See also the 1991 work of Chan and Mathew \cite{Chan2}.

\bigskip

The related idea of computing a low-rank approximation or ``skeleton" of a matrix by means of randomized testing, albeit without a priori knowledge of the row and column spaces, was extensively studied in recent work of Rokhlin et al. \cite{randomSVD1, randomSVD3}, and Martinsson and Tropp \cite{randomSVD2}.

\section{Numerical results}\label{sec:num}

The classical Marmousi benchmark example is the basis of all our numerical experiments. The forward model is taken to be the linearized wave equation
\[
m_0(x) u_{tt} - \Delta u = - \delta m(x) (u_0)_{tt},
\]
where the incident field $u_0$ obeys
\[
m_0(x) (u_0)_{tt} - \Delta u_0 = f_s,
\]
with $f_s(x,t) = \delta(x-s) w(t)$. The wavelet $w(t)$ is taken to be the second derivative of a gaussian (Ricker wavelet). The background medium $m_0$ is either taken to be constant (in Sections \ref{sec:numbasic}, \ref{sec:numgen}), or a smoothed version of the original Marmousi model with various degrees of smoothing (in Section \ref{sec:numvariable}). The data $d(r,s,t)$ are then collected as the samples of $u$ at receiver positions $r$ and source positions $s$ at the surface $z=0$, and all adequate times $t$.

\bigskip

The same equations are then used for the imaging, with $m_0$ and $f_s$ assumed known, but not $\delta m(x)$. This is known as the ``inversion crime", as any real-life imaging application would also require to solve for $m_0$ and $f_s$ -- problems that we leave aside in this paper. Notice also that the forward model is \emph{linear} in $\delta m(x)$, a clearly uncalled-for assumption in practice since it neglects multiple scattering. A better wave equation for $u$ would have $u_{tt}$ in place of $(u_0)_{tt}$ in the right-hand-side. We nevertheless made this assumption so as not to obscure the fact that the Hessian is intrinsically present to correct the solution of the linearized inverse problem.

\bigskip

For the convenience of being able to run hundreds of simulations in a matter of hours, we choose to consider a 2D problem on a square domain with $N^2$ points, $N = 127$ for most of the results shown. A perfectly matched layer (PML) of width $.15 N$ surrounds the domain of interest. The numerical method has spectral differences in space, and second-order differences in time. The poststack imaging operator $F^*$ performs a stack on three sources maximally spaced from each other (albeit not in the PML). More sources were used in some of the numerical experiments, but this did not significantly affect the inverse Hessian. As is well-known, the main advantage of using more sources is the robustness to noise. (All the imaging results are robust to additive gaussian white noise, but not to purely multiplicative gaussian white noise.)

\bigskip

Two types of preconditioners are compared:
\begin{itemize}
\item {\bf Rn}: Fitting of the inverse Hessian from randomized curvelet trial functions. This preconditioner is denoted as Rn where $n$ is the number of trial functions used for the fitting, e.g. R4 is four functions were used.

\item {\bf Kn}: Fitting the inverse Hessian from trial functions taken in the Krylov subspace of the migrated image. This preconditioner is denoted Kn where $n$ is the number of trial functions used for the fitting, e.g. K2 if both the migrated image and the remigrated image were used. This is essentially the method of Nammour and Symes \cite{Nammour1, Nammour2}, with the slight improvement of using the full expansion (\ref{eq:compr2}) in place of (\ref{eq:compr1}) -- a minor point.

\end{itemize}
\bigskip

In both cases the $B_i$ are the elementary symbols of equation (\ref{eq:compr2}). Different numbers of terms are tested in this pseudodifferential expansion: in order of decreasing importance, the parameters are 1) number of Fourier modes in $x$, 2) number of Fourier modes in the wavevector argument $\theta$, and 3) number of Chebyshev modes in the wavenumber $|k|$. The right balance of parameters in each dimension was obtained manually for best accuracy; only their total number (their product) is reported.

\bigskip

The action of the preconditioners on the migrated image $F^* d$ is compared to the image obtained after 200 gradient descent steps for the (linearized) least-squares functional. The refinement of this brute force method to an iterative solver such as GMRES or LSQR is important in practice, but was not investigated in the scope of this paper.

\bigskip

Errors between models are measured in the relative mean-squared sense, i.e. if $\delta m_1$ is a reference model and $\delta m_2$ another model, then
\[
\mbox{MSE}(\delta m_1,\delta m_2) = \frac{\| \delta m_1 - \delta m_2 \|_2}{\| \delta m_1 \|_2}.
\]

\subsection{Basic results}\label{sec:numbasic}

The action of the preconditioners on the migrated image is satisfactory: as the figures below show it is visually closer to the image obtained after 200 gradient steps than the migrated image. 

\begin{figure}[H]
\centering
\includegraphics[width=5.2cm]{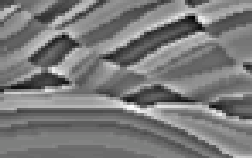}\hfill
\includegraphics[width=5.2cm]{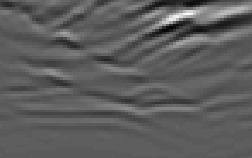} \hfill
\includegraphics[width=5.2cm]{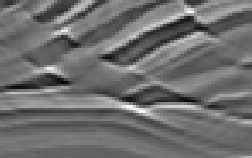} \hfill
\caption{Left: oscillatory wave speed profile (``reflectors") used to produce wavefield data. The forward model is the linearized wave equation with a unit background speed. Middle: migrated image, obtained by reverse-time migration. Right: image obtained by 200 gradient descent steps to solve the linearized least-squares problem.}
\label{fig:num1}
\end{figure}

\begin{figure}[H]
\centering
\includegraphics[width=5.2cm]{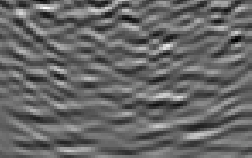}\hfill
\includegraphics[width=5.2cm]{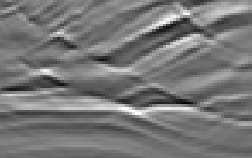}\hfill
\includegraphics[width=5.2cm]{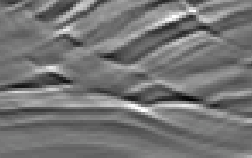}\hfill
\caption{Left: a randomized curvelet trial function, used for testing the Hessian in order to fit the inverse Hessian. Middle: image obtained by applying the R4 preconditioner to the migrated image. Right: image obtained by applying the K1 preconditioner to the migrated image.}
\end{figure}

The Krylov preconditioner K1 usually works well on the migrated image. The randomized preconditioner R1 is often a notch worse than K1, but when going up to R4 and higher the performance becomes very comparable to K1. We did not find an instance where any Rn, regardless of $n$, would significantly outperform K1 (a puzzling observation). However, we notice in Figure \ref{fig:num4} that as the dimension of the Krylov subspace increases, the performance of K2, K3, etc. deteriorates very quickly. This is in contrast to what was advocated in \cite{Nammour1, Nammour2}. 

\bigskip

\begin{figure}[H]
\centering
\includegraphics[width=8cm]{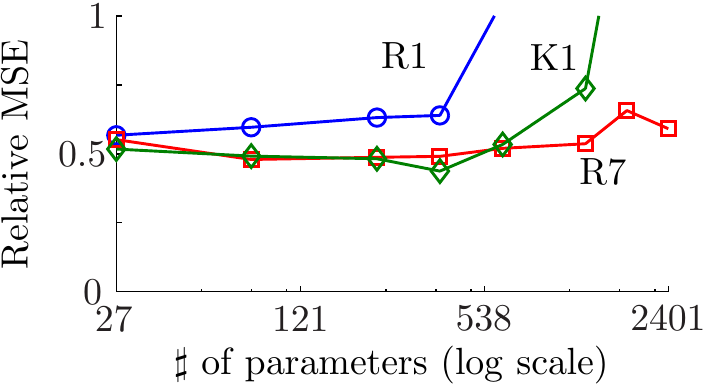}\hfill
\includegraphics[width=8cm]{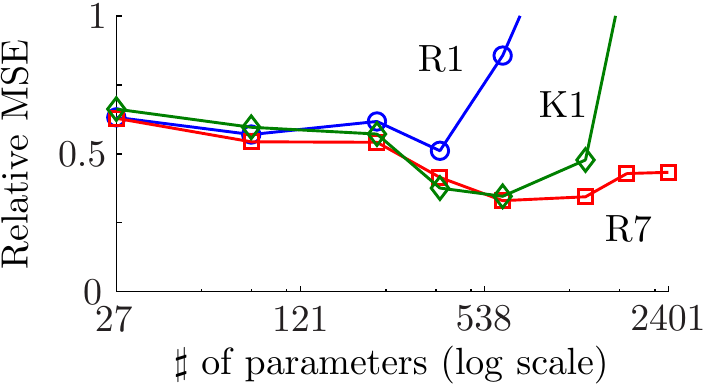}
\caption{Relative MSE of the R1, R7 and K1 preconditioners, as a function of the number of parameters. Left: preconditioned migrated image vs. slightly modified ``true" image of Figure \ref{fig:num1}, left. By ``slightly modified", we mean that a curvelet mask is taken to only measure the components of those images in the set $S$ (see section \ref{sec:curvelettrial}). Right: preconditioned migrated image vs. recovered image of Figure \ref{fig:num1}, right. No curvelet mask is taken here. Any MSE below 1 (100 percent relative error) indicates that the preconditioning is working.}
\label{fig:num3}
\end{figure}

There is a sweet spot in the number of parameters in the symbol expansion of the inverse Hessian, around 500 to 1000 for the numerical scenario considered. See Figure \ref{fig:num3}. If the number of parameters is too small, the inverse Hessian is not properly represented. If the number of parameters is too large, they are either not used to improve the representation of the Hessian, or their large number leads to ill-conditioning of the fitting problem (hence large numerical errors.)

Note that in this experiment the Hessian is a $16,384$-by-$16,384$ matrix. Its numerical rank hovers in the few thousands; more precisely, for a top singular value normalized to unity, the $\eps$-rank as a function of $\eps$ is given by the following table. We attribute the rank deficiency mostly to the perfectly matched layer (PML) and other windows applied.

\begin{center}
\begin{tabular}{c|c}
$\eps$ & $\eps$-rank \\ \hline
1e-1 & 435 \\
1e-2 & 1367 \\
1e-3 & 2164 \\
1e-4 & 2803 \\
1e-5 & 3250 \\
1e-6 & 3624
\end{tabular}
\end{center}

\subsection{Generalization error}\label{sec:numgen}

\begin{figure}[H]
\centering
\includegraphics[width=8cm]{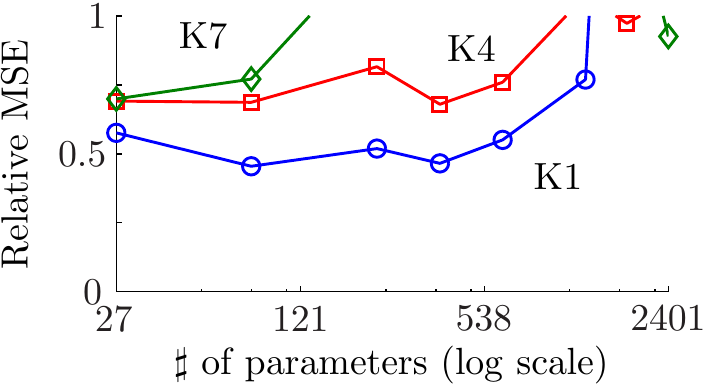}\hfill
\includegraphics[width=8cm]{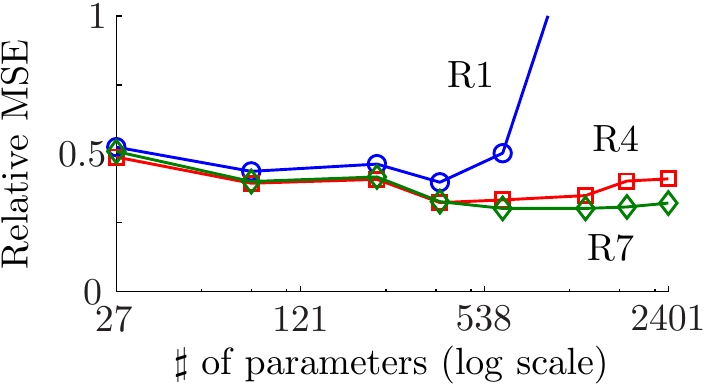} \hfill
\caption{Relative MSE of generalization. The setup is the same as in the previous section, except that the Marmousi model was replaced by an (indepently drawn) randomized curvelet trial function. Here the error is simply the relative MSE for the reconstruction of this randomized trial function, from applying the Hessian followed by a preconditioner. The x axis shows the number of parameters. Left: K1, K4, and K7 preconditioners applied to the randomized trial function, vs. image obtained by 200 gradient descent iterations. The performance quickly degrades with the order of the Krylov subspace. Right: R1, R4, and R7 preconditioners applied to the randomized curvelet trial function, vs. reference image. Notice that the error is smaller than in the Krylov case, and that the performance does not degrade with the index of the preconditioner.}
\label{fig:num4}
\end{figure}

The Rn preconditioners show their true potential when the inverse Hessian is applied to another randomized trial function, drawn independently from those used for fitting the symbol, see Figure \ref{fig:num4}, right. Generalizability to a large linear subspace of models is as the theory predicts. The Krylov preconditioners, on the other hand, show some fragility here. They are not designed to work when applied on images far from the remigrated image, and indeed, the error level is higher for K1 than for any Rn.

The degradation of the Kn preconditioners as $n$ increases is understandable. In applying the normal operator $1, 2, \ldots, n$ times to the migrated image, information is lost in all but the eigenspaces corresponding to leading eigenvalues. This is well-known from the analysis of the power method in linear algebra. As a result, the disproportionate weight lent to those subspaces ``hijacks" most of the degrees of freedom of the symbol expansion and prevents a good fit.

\bigskip

The robustness of the Rn preconditioners offered by generalizability may be useful in the scope of preconditioned gradient descent iterations. While $H^{-1}$ is applied to $F^*d$ (migrated image) in the first iteration, it is subsequently applied to $F^*(d - F \, \delta  m_k)$ (migrated residual). The latter will deviate from $F^* d$ in the course of the iterations, resulting in a weaker K1 preconditioner.


\subsection{Variable media}\label{sec:numvariable}

The curvelet mask used in the definition of the randomized trial functions is a set $S$ in curvelet space indicating whether the curvelet is ``visible in the dataset" or not. In the case of a uniform medium, this information is obtained by considering the fan of couples of lines originating from each curvelet's center point, for which the angle of incidence equals the angle of reflection. For a given curvelet the test is whether one of the lines joins the curvelet to a source while the other line joins the curvelet to a receiver. If this test returns a positive match for one couple of lines, we declare that the curvelet is active and its index belongs to the set $S$.

\bigskip

In the case of smooth variable media, the test is similar but now involves ray tracing, i.e., computing the trajectories of the Hamiltonian system of geometrical optics. This is performed ray-by-ray using the high-order adaptive Runge-Kutta time integrator ode45 built in Matlab. Ray-tracing is normally not a computational bottleneck; if solving for the rays one-by-one is too slow, a fast algorithm such as the phase-flow method of Ying and Cand\`{e}s \cite{PFM} can be set up to speed up the process.

\bigskip

For the numerical experiment we take the smooth part of the Marmousi model $M(x)$ and smooth it further by convolution with a radial bump. This operation is realized in the wavevector domain, by multiplying the Fourier transform of $M(x)$ by the indicator function of a disk of radius $rN$ (the whole wavevector space is a square of sidelength $N$). We let $0 \leq \gamma \leq 0.4$ and consider $M_\gamma(x)$ the further-smoothed Marmousi background model velocity. Then we set
\[
m_0(x) = \left( 1 - \frac{\gamma}{0.4} \right) \int M(x) dx \, + \, \frac{\gamma}{0.4} M_\gamma(x).
\]
If $\gamma = 0$ we recover a uniform medium. The MSE of the R5 preconditioner as a function of $0 \leq \gamma \leq 0.4$ is shown below. Most of the numerical tests performed in the earlier sections were repeated in variable media: we did not find that any particular plot was worth reporting, as the performance systematically degrades in a predictable manner as $\gamma$ increases.

\begin{figure}[H]
\centering
\includegraphics[width=8cm]{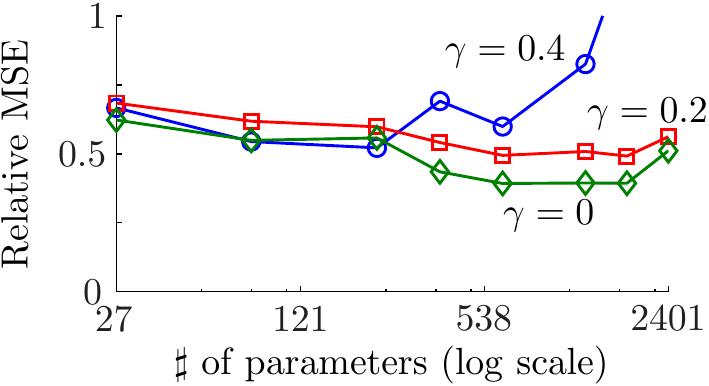} 
\caption{Relative MSE of the R5 preconditioner applied to the migrated image, vs. the image obtained by 200 gradient descent iterations (Figure \ref{fig:num1}, right). The x axis shows the number of parameters. The different curves refer to different smoothness levels of the model velocity, as explained in the text.}
\end{figure}

\subsection{Other tests}

Other sizes, from $N = 64$ to $N = 256$ were tested and showed similar performance levels. 

Other randomized trial functions than ``curvelet-masked noise" were attempted, such as
\begin{itemize}
\item Gaussian white noise in model space, which failed badly because it contains too much energy in the nullspace, with high probability.
\item Gaussian white noise in data space, migrated to model space. Such trial functions still have too much energy in the nullspace and led to unequivocally poor results.
\item Gaussian white noise in model space, to which the normal operator is applied. These trial functions work well, and show error levels comparable (at times slightly worse) than the curvelet trial functions. They have the advantage of being simple to define -- no need for curvelets -- but more complicated to compute as each randomized trial function requires one application of the expensive Hessian.
\item Gaussian white noise in model space, to which the normal operator is applied, followed by a diagonal operation in curvelet space where the coefficient magnitudes are either put to 1 or to zero if they are under a small threshold. Coefficient phases are unchanged. These trial functions are comparable to the simpler ones defined directly in curvelet space.
\item Other distributions than gaussian for the noise: this did not give rise to any noticeable difference in our numerical experiments. Lemmas are indeed often available to pass from one distribution to the other in large deviation theory.
\end{itemize}

The fitting of the inverse Hessian was also realized from an application of the Hessian to the desired unknown model that served to create the data. This operation can of course not be performed in practice since we are precisely trying to invert for this unknown model. But the numerical experiment is very instructive: it shows that the relative MSE of the Rn preconditioner applied to the migrated image decays to such small values as 0.1 when the number of parameters is large enough; the MSE does not stall on a plateau at 0.3 like it does in all the figures above. This goes to show that the pseudodifferential expansion is instrinsically good, but that neither the Krylov fit nor the randomized fit is fine enough to predict the right coefficients. This leaves exciting room for improvement of the method.

\section{Theory}\label{sec:theory}

\subsection{Invertibility of $M$}
To carry out the least squares minimization in Section \ref{sec:randfit}, the $n$ by $n$ matrix $M$ has to be well-conditioned. In this section, we will show that this happens with high probability (whp) when the number of parameters $p$ is related to the (numerical) rank $r$ of $H$ through 
\[
r \geq C p^2 \log p, \qquad \mbox{ for some } C > 0.
\]

If $H$ were an invertible matrix, we would simply let $y \sim N(0,1)^n$, independent and identically distributed (iid). But in the general case, and as mentioned earlier, we should make sure that $y$ is properly ``colored" to avoid the nullspace of $H$. While our numerical solution to this problem is approximate, we will assume for simplicity
that we can exactly project $y$ onto the range space of $H$,
\[
\tilde{\mathbf{y}} = P \mathbf{y},
\]
where $P$ is the orthogonal projector onto Ran$(H)$. 

The random matrix $M$ to invert for the fitting step is then
\[
M_{ij} = \tilde{\mathbf{y}}^T (H B_i^T B_j H) \tilde{\mathbf{y}}.
\]
It holds that $M_{ij}= \mathbf{y}^T (H B_i^T B_j H) \mathbf{y}$ without the tildes, hence
\[
\E M_{ij} = \mbox{Tr}(HB_i^T B_j H).
\]
It is assumed that $\E M$ is positive definite and well-conditioned; our argument consists in showing that $M$ does not depart too much from its expectation whp.

Let $\|\cdot\|$ and $\|\cdot\|_F$ denote the spectral and Frobenius norms respectively. We denote by $\kappa$ the condition number of $\E M$,
\[
\kappa = \| \E M \| \, \| (\E M)^{-1} \|.
\]
We also need to consider $\eta > 0$, the smallest number such that
\[
\| HB_i \| \leq \frac{\eta}{\sqrt{r}} \| HB_i \|_F,
\]
uniformly over $i$. We may call $\eta$ the ``weak condition number" of the collection of $H B_i$.

Both $\kappa$ and $\eta$ are greater than 1, but it will be manifest from the way they enter the estimates below that they ought to be small (close to 1). If $\eta$ is small, then $HB_i$ has approximate numerical rank $r$, i.e., the largest $r$ singular values are comparable in size.



The following result is a perturbative analysis quantifying the size of $\| M-\E M \|$ in relation to $\| \E M \|$. 

\begin{theorem}\label{thm:psquare} Assume that $H$ is a symmetric rank-$r$ matrix that can be written as $H = \sum_{i = 1}^p c_i B_i$. Define $M_{ij}$ and $\eta$ as above. For all $0 < \eps \leq 1$, there exists a number $C(\eps, \eta) > 0$ such that, if
\[
r \geq C(\eps, \eta) \, p^2 \log p,
\]
then with high probability
\[
\| M - \E M \| < \eps \| \E M \|.
\]
Explicitly, $C(\eps, \eta) = 160 \eta^4 \eps^{-2}$, and the ``high probability"  is at least $1 - 2 p^{-8}$.
\end{theorem}

Before we prove this theorem, let us explain how invertibility of $M$ follows at once. Since the condition number of $\E M$ is $\kappa$, its minimum eigenvalue obeys
\[
\lambda_{\min}(\E M) \geq \frac{1}{\kappa} \, \| \E M \|.
\]
When a matrix is perturbed, the change in eigenvalues is controlled by the spectral norm of the perturbation, so
\[
\lambda_{\min}(M) \geq \left( \frac{1}{\kappa} - \eps \right) \, \| \E M \|, \qquad \mbox{whp}.
\]
It suffices therefore to apply the theorem above with $\eps < \frac{1}{\kappa}$ to ensure invertibility of $M$.




\begin{proof}[Proof of Theorem \ref{thm:psquare}]

Let us first settle that $M_{ij}= \mathbf{y}^T (H B_i^T B_j H) \mathbf{y}$, without the tildes. It suffices to argue that $HP = H$. By transposition, and symmetry of both $H$ and $P$, it suffices to show that $PH = H$. This latter equation is obviously true since $P$ acts as the identity on the range space of $H$.

\bigskip

Now let $L=\|\E M\|$. Our proof considers the statistics of $M_{ij}$ element-wise as a quadratic form of the gaussian random vector $\mathbf{y}$. We will show that $M_{ij}$ is highly unlikely to be more than $\eps L/p$ away from $\E M_{ij}$. In what follows we use the $\|\cdot\|_1$ and $\|\cdot \|_{\infty}$ induced matrix norms -- the maximum absolute column and row sums respectively. If we can show that $|M_{ij}-\E M_{ij} |<\eps L/p$ for all $i,j$, then the following inequality completes the proof:
$$\|M-\E M\|_2 \leq \frac{1}{2} \left( \|M-\E M\|_1 + \|M-\E M\|_{\infty}\right)\leq \eps L.$$

The statistics of quadratic forms $\mathbf{y}^T A \mathbf{y}$ were perhaps first completely studied by Grenander, Pollak and Slepian \cite{Grenander59}. In a nutshell , the variance of the quadratic form $M_{ij} = \mathbf{y}^T (H B_i^T B_j H) \mathbf{y}$ is known to be proportional to $\|\frac{1}{2}(H B_i^T B_j H + H B_j^T B_i H)\|_F^2$. We seek to bound these variances using the fact that the $HB_i$ are ``weakly well-conditioned".

Fix $i$. We know that
$$L =\| \E M\| \geq |\E M_{ii}| = |\tr(H B_i^T B_i H)| = \|H B_i\|_F^2.$$

Using the definition of $\eta$ we obtain a stronger bound on the spectral norm, namely $\|H B_i\|\leq \eta \|HB_i\|_F/\sqrt{r} \leq \eta \sqrt{L/r}.$ The implication is that for all $i,j$,
\begin{equation}\label{eq:var2}
\| H B_i^T B_j H\|_2 \leq \eta^2 L/r .
\end{equation}
As for the Frobenius norm of $H B_i^T B_j H$, we make use of the fact that $H$ has rank $r$ to bound
\begin{equation}\label{eq:varF}
\|H B_i^T B_j H\|_F \leq \|H B_i^T B_j H\|\sqrt{r} = \eta^2 L/\sqrt{r}.
\end{equation}

We are now ready to bound $\Pr \left( |M_{ij}-\E M_{ij}|>\eps L(p)/p \, \right)$. For clarity, fix $i,j$ and let $A=H B_i^T B_j H$.  The standard deviation of $A$ is proportional to $\|A\|_F$, which by Eq. (\ref{eq:varF}) is roughly on the order of $L/\sqrt{r}$ or $L/p$. This is qualitatively correct. For an explicit bound, we refer to Bechar \cite{Bechar09}, who builds on the work of \cite{Grenander59} to state the following.

\begin{lemma}\label{lemma:bechar}
Let $A\in \reals^{n \times n}$ and $y\sim N(0,1)^n$ iid. Then for any $\lambda>0$,
$$\Pr(|\mathbf{y}^T A \mathbf{y} - \E \mathbf{y}^T A \mathbf{y}| \geq \|A+A^T\|_F \sqrt{\lambda} + 2\|A\| \lambda) \leq 2 \exp(-\lambda).$$
\end{lemma}


We pick $\lambda = 10 \log p$. It is straightforward to verify that with this choice of $\lambda$, with the definition of $C(\eps, \eta)$, and with equations (\ref{eq:var2}) and (\ref{eq:varF}), we have
\[
\|A+A^T\|_F \sqrt{\lambda} + 2 \|A\|_2 \lambda \leq \eps L/p.
\]
It follows that 
$$\Pr(|M_{ij}-\E M_{ij}|\geq \eps L(p)/p) \leq 2 \exp(-\lambda) < 2 p^{-10}.$$

An union bound over $p^2$ pairs of $i,j$'s concludes the proof. Note in passing that we made no effort to minimize $C(\eps, \eta)$.

\end{proof}

Finally, we sketch a standard procedure to handle complex-valued matrices. Instead of taking the symmetric part of $A$ by $\mathbf{y}^T A \mathbf{y}=\mathbf{y}^T (\frac{1}{2}(A+A^T))\mathbf{y}$, decompose it into Hermitian and anti-Hermitian components, that is $\mathbf{y}^T A \mathbf{y} = \mathbf{y}^T A_1 \mathbf{y} -i \mathbf{y}^T A_2 \mathbf{y}$ where $A_1=\frac{1}{2}(A+A^*)$ and $A_2=\frac{i}{2}(A-A^*)$ are both Hermitian. Then bound the deviations from their expectations separately by $|\mathbf{y}^T A \mathbf{y}-\E \mathbf{y}^T A \mathbf{y}| \leq |\mathbf{y}^T A_1 \mathbf{y} - \E \mathbf{y}^T A_1 \mathbf{y}|+|\mathbf{y}^T A_2 \mathbf{y} - \E \mathbf{y}^T A_2 \mathbf{y}|$. Repeat similar arguments and invoke Lemma \ref{lemma:bechar} to show that each term is less than $\eps L/2p$ whp.

\subsection{Rationale for curvelets}\label{sec:theory2}

The success of the proposed method for inverting the Hessian depends on the property of phase-space localization of curvelets. Good localization of a basis function like a curvelet near a point $(x,k)$ implies that it will only ``see" values of the symbol $a(x,k)$ near that point, when acted upon by the Hessian.

The following result makes this heuristic precise; it is a minor modification of a theorem of Stolk \cite{HerrmannStolk} so the proof is omitted.

\begin{theorem} (Stolk, 2008).
Let $a(x,k)$ be the pseudodifferential symbol of the wave equation Hessian $H$, as in equation (\ref{eq:symbol}), and assume that it obeys (\ref{eq:type10}) with $m=1$. Consider the zeroth-order symbol $a(x,k) |k|^{-1}$ of the operator $H (- \Delta)^{-1/2}$. Denote by $\tilde{H} (- \Delta)^{-1/2}$ the diagonal approximation of $H (- \Delta)^{-1/2}$ in curvelet space, with the sampled symbol as multiplier,
\[
\tilde{H} (- \Delta)^{-1/2} f = \sum_{\mu} \vf_\mu(x) a(x_\mu, k_\mu) |k_\mu|^{-1} \int \overline{\vf_\mu}(x) f(x) \, dx.
\]
If $f$ obeys $\hat{f}(k) = 0$ for $|k| \leq k_{\min}$, then there exists $C > 0$ such that
\[
\| ( \tilde{H} - H ) (- \Delta)^{-1/2} f \|_2 \leq \frac{C}{\sqrt{k_{\min}}} \| f \|_2.
\]
\end{theorem}

In other words, the more oscillatory the model $f(x)$ the better the diagonal approximation of the Hessian via curvelets. Hence the larger $k$ the better the ``probing" character of a curvelet near its center in phase-space.

\bigskip

The theorem above is also true for another frame of functions, the wave atoms of Demanet and Ying \cite{waveatoms}, but would not be true for wavelets, directional wavelets, Gabor functions, or ridgelets.

\section{Conclusion}

This paper presents a preconditioner for the wave equation Hessian based on ideas of randomized testing, pseudodifferential symbols, and phase-space localization. Numerical experiments show that the proposed solution belongs to a class of effective ``probing" preconditioners. The precomputation only requires applying the wave equation Hessian once, or a small number of times.

Fitting the inverse Hessian involves solving a small least-squares problem, of size $p$-by-$p$, where $p$ is much smaller than $n$ and the Hessian is $n$-by-$n$. Even if $p$ were on the order of $n$ the proposed method would be very advantageous since constructing each row of the Hessian requires going back to the much higher dimensional data space.

It is anticipated that the techniques developed in this paper will be of particular interest in 3D seismic imaging and with more sophisticated physical models that require identifying a few different parameters (elastic moduli, density). In that setting, properly inverting the Hessian with low complexity algorithms to unscramble the multiple parameters will be particularly desirable.



\end{document}